\documentclass[11pt,leqno,fleqn]{article}
\usepackage{epsf,amsmath,amssymb,amscd,hyperref,srcltx}
\usepackage{amsfonts}
\usepackage{psfrag}
\newcommand{\dps}{\displaystyle}
\newcommand{\T}{^{\sf T}}
\newcommand{\dy}[2]{%
\refstepcounter{equation}%
\label{#1}%
\begin{list}{}{
\topsep 5mm
\leftmargin 18mm
\rightmargin 0cm
\itemsep 0mm
\listparindent 0mm
\parsep 0mm
\itemsep 0mm
\labelsep 0mm
\labelwidth 18mm
}%
\item[\rm (\theequation)\hfill]
#2
\end{list}%
}
\newcommand{\dyy}[2]{\dy{#1}{\raggedright$\dps#2$}}
\newcommand{\de}[2]{\dy{#1}{\raggedright$\displaystyle #2 $}}
\newcounter{hulpstelling}
\newcommand{\lemma}[2]{\refstepcounter{hulpstelling}\vspace{4mm}\noindent{\bf Lemma \thehulpstelling.}\label{#1}{\it #2}}
\newcommand{\pf}{\vspace{3mm}\noindent{\bf Proof.}\ }
\newcommand{\bx}{\hspace*{\fill} \hbox{\hskip 1pt \vrule width 4pt height 8pt depth 1.5pt \hskip 1pt}

\addvspace{4mm}}
\newcommand{\bxx}{\hspace*{\fill} \hbox{\hskip 1pt \vrule width 4pt height 8pt depth 1.5pt \hskip 1pt}}
\newcommand{\rf}[1]{{\rm (\ref{#1})}}
\newcommand{\oN}{{\mathbb{N}}}
\newcommand{\oR}{{\mathbb{R}}}
\newcommand{\Tr}{{\text{Tr}}}
\begin{document}

\begin{center}
{\LARGE\bf A Pythagoras proof of Szemer\'edi's regularity lemma

}

\medskip
{\large
Notes for our seminar  --- Alexander Schrijver \footnote{ CWI and University of Amsterdam.
Mailing address: CWI, Science Park 123, 1098 XG Amsterdam,
The Netherlands.
Email: lex@cwi.nl.}}

\end{center}

\noindent
{\small{\bf Abstract.}
We give a short proof of Szemer\'edi's regularity lemma, based on
elementary geometry.

}

\medskip
The `regularity lemma' of Endre Szemer\'edi [1] roughly asserts that, for each
$\varepsilon>0$, there exists a number $k$ such that the vertex set $V$ of
any graph $G=(V,E)$ can be partitioned into at most $k$ {\em almost} equal-sized classes so
that between {\em almost} any two classes, the edges are distributed {\em almost}
homogeneously.
Here {\em almost} depends on $\varepsilon$.
The important issue is that $k$ (though generally extremely huge) only depends on $\varepsilon$,
and not on the size of the graph.
The lemma has several applications in graph and number theory, discrete geometry, and
theoretical computer science.

We give a short proof based on elementary Euclidean geometry.
The general line of the proof is like that of the standard proof (in fact, Szemer\'edi's original proof), but
most of the technicalities are swallowed by Pythagoras' theorem.
We prove two lemmas, one on `$\varepsilon$-balanced' partitions,
the other on `$\varepsilon$-regular' partitions.

Let $V$ be a finite set.
A {\em partition} of $V$ is a collection of disjoint nonempty sets (called {\em classes}) with union $V$.
Partition $Q$ of $V$ is a {\em refinement} of partition $P$ if each class of $Q$ is contained
in some class of $P$.
For $\varepsilon>0$, partition $P$ of $V$ is called {\em $\varepsilon$-balanced}
if $P$ contains a subcollection $ C$
such that all sets in $C$ have the same size and such that
$|V\setminus\bigcup C|\leq\varepsilon|V|$.

\lemma{20no12h}{
Each partition $P$ of $V$ has an $\varepsilon$-balanced refinement $Q$
with $|Q|\leq (1+\varepsilon^{-1})|P|$.
}

\pf
Define $t:=\varepsilon|V|/|P|$.
Split each class of $P$ into classes, each of size $\lceil t\rceil$,
except for at most one of size less than $t$.
This gives $Q$.
Then $|Q|\leq |P|+|V|/t=(1+\varepsilon^{-1})|P|$.
Also, the union of the classes of $Q$ of size less than $t$ has
size at most $|P|t=\varepsilon|V|$.
So $Q$ is $\varepsilon$-balanced.
\bx

Let $G=(V,E)$ be a graph.
For nonempty $I,J\subseteq V$, the {\em density}
$d(I,J)$ of $(I,J)$ is the number of adjacent pairs of
vertices in $I\times J$, divided by $|I\times J|$.
Call the pair $(I,J)$ {\em $\varepsilon$-regular} if for all $X\subseteq I, Y\subseteq J$:
\dy{27no12b}{
if $|X|>\varepsilon|I|$ and $|Y|>\varepsilon|J|$ then $|d(X,Y)-d(I,J)|\leq\varepsilon$.
}
A partition $P$ of $V$ is called {\em $\varepsilon$-regular} if
\de{20no12f}{
\sum_{I,J\in P\atop (I,J)\text{ \rm $\varepsilon$-irregular}}
\hspace*{-4mm}
|I||J|\leq\varepsilon |V|^2.
}

\vspace*{-2mm}
For Lemma \ref{28no12a} we need the following.
Consider the matrix space $\oR^{V\times V}$, with the Frobenius norm
$\|M\|=\Tr(M\T M)^{1/2}$ for $M\in\oR^{V\times V}$.
For nonempty $I,J\subseteq V$,
let $L_{I,J}$ be the $1$-dimensional subspace of $\oR^{V\times V}$ 
consisting of all matrices that are constant on $I\times J$ and 0
outside $I\times J$.
For any $M\in\oR^{V\times V}$, let $M_{I,J}$ be the orthogonal projection of $M$
onto $L_{I,J}$.
So the entries of $M_{I,J}$ on $I\times J$ are all equal
to the average value of $M$ on $I\times J$.

If $P$ is a partition of $V$, let $L_P$ be the sum
of the spaces $L_{I,J}$ with $I,J\in P$,
and let $M_P$ be the orthogonal projection of $M$ onto $L_P$.
So $M_P=\sum_{I,J\in P}M_{I,J}$.
Note that if $Q$ is a refinement of $P$, then $L_P\subseteq L_Q$, hence
$\|M_P\|\leq\|M_Q\|$.

\lemma{28no12a}{
Let $\varepsilon>0$ and $G=(V,E)$ be a graph,
with adjacency matrix $A$.
Then each $\varepsilon$-irregular partition $P$ has
a refinement $Q$
with $|Q|\leq |P|4^{|P|}$ and
$\|A_Q\|^2>\|A_P\|^2+\varepsilon^5|V|^2$.
}

\pf
Let $(I_1,J_1),\ldots,(I_n,J_n)$ be the $\varepsilon$-irregular pairs in $P^2$.
For each $i=1,\ldots,n$,
we can choose (by definition \rf{27no12b}) subsets $X_i\subseteq I_i$
and $Y_i\subseteq J_i$ with $|X_i|>\varepsilon|I_i|$,
$|Y_i|>\varepsilon|J_i|$ and $|d(X_i,Y_i)-d(I_i,J_i)|>\varepsilon$.
For any fixed $K\in P$, there exists a partition $Q_K$ of $K$ such that each $X_i$
with $I_i=K$ and each $Y_i$ with $J_i=K$ is a union of classes of $Q_K$ and such that
$|Q_K|\leq 2^{2|P|}=4^{|P|}$.
\footnote{For any collection $C$ of subsets of a finite set $S$, there is a
partition $R$ of $S$ such that any set in $C$ is a union of classes of
$R$ and such that $|R|\leq 2^{|C|}$:
take $R:=
\{\bigcap_{X\in D}X\cap\bigcap_{Y\in C\setminus D}S\setminus Y\mid D\subseteq C\}\setminus\{\emptyset\}$.}
Let $Q:=\bigcup_{K\in P}Q_K$.
Then $Q$ is a refinement of $P$ such that each $X_i$ and each
$Y_i$ is a union of classes of $Q$.
Moreover, $|Q|\leq |P|4^{|P|}$.

Now note that for each $i$, since
$(A_Q)_{X_i,Y_i}=A_{X_i,Y_i}$ (as $L_{X_i,Y_i}\subseteq L_Q$) and since
$A_{X_i,Y_i}$ and $A_P$ are constant on
$X_i\times Y_i$, with values $d(X_i,Y_i)$
and $d(I_i,J_i)$, respectively:
\de{14de12a}{
\hspace*{-4mm}
\|(A_Q\hspace*{-0.6mm}-\hspace*{-0.7mm}A_P)_{X_i,Y_i}\|^2 =
\|A_{X_i,Y_i}\hspace*{-0.6mm}-\hspace*{-0.6mm}(A_P)_{X_i,Y_i}\|^2
=
|X_i||Y_i|(d(X_i,Y_i)\hspace*{-0.6mm}-\hspace*{-0.6mm}d(I_i,J_i))^2
>
\varepsilon^4|I_i||J_i|.
}
Then negating \rf{20no12f} gives with Pythagoras, as
$A_P$ is orthogonal to $A_Q-A_P$ (as $L_P\subseteq L_Q$),
and as the spaces
$L_{X_i,Y_i}$ are pairwise orthogonal,
\dyy{25no12b}{
\hspace*{-4mm}
\|A_Q\|^2\hspace*{-0.7mm}-\hspace*{-0.7mm}\|A_P\|^2
=
\|A_Q\hspace*{-0.6mm}-\hspace*{-0.7mm}A_P\|^2
\geq
\sum_{i=1}^n
\|(A_Q\hspace*{-0.6mm}-\hspace*{-0.7mm}A_P)_{X_i,Y_i}\|^2
\geq
\sum_{i=1}^n
\varepsilon^4|I_i||J_i|
>
\varepsilon^5|V|^2.
\bxx
}

Define $f_{\varepsilon}(x):=(1+\varepsilon^{-1})x4^x$ for $\varepsilon,x>0$.
For $n\in\oN$, $f_{\varepsilon}^n$ denotes the $n$-th iterate of $f_{\varepsilon}$.

\vspace{4mm}
\noindent
{\bf Szemer\'edi's regularity lemma.}{
For each $\varepsilon>0$ and graph $G=(V,E)$, each partition $P$ of
$V$ has an $\varepsilon$-balanced $\varepsilon$-regular refinement of size
$\leq f_{\varepsilon}^{\lfloor\varepsilon^{-5}\rfloor}((1+\varepsilon^{-1})|P|)$.
}

\pf
Let $A$ be the adjacency matrix of $G$.
Starting with $P$, iteratively apply Lemmas \ref{20no12h} and \ref{28no12a} alternatingly.
At each application of Lemma \ref{20no12h}, $\|A_P\|^2$ does not decrease,
and at each application of Lemma \ref{28no12a}, $\|A_P\|^2$ increases by more that
$\varepsilon^5|V|^2$.
Now, for any partition $Q$ of $V$, $\|A_Q\|^2\leq \|A\|^2\leq |V|^2$.
Hence,
after at most $\lfloor\varepsilon^{-5}\rfloor$ iterations we must have an
$\varepsilon$-balanced $\varepsilon$-regular partition as required.
\bx

We note that if
$P$ is an $\varepsilon$-balanced $\varepsilon$-regular partition of $V$,
and $C\subseteq P$ is such that all sets in $C$ have
the same size and such that
$|V\setminus\bigcup C|\leq\varepsilon|V|$, then the number $s$ of
$\varepsilon$-irregular pairs in $C^2$ is at most
$\varepsilon(1-\varepsilon)^{-2}|C|^2$.
For let $t$ be the common size of the sets in $C$.
Then, by \rf{20no12f},
$st^2\leq\varepsilon|V|^2\leq
\varepsilon(1-\varepsilon)^{-2}|\bigcup C|^2
=
\varepsilon(1-\varepsilon)^{-2}(t|C|)^2
=
\varepsilon(1-\varepsilon)^{-2}|C|^2t^2
$.

\bigskip
\noindent
{\bf Reference}

\medskip
\noindent
\small
[1] E. Szemer\'edi, 
Regular  partitions of graphs,
in: {\em Probl\`emes combinatoires et th\'eorie des graphes}
(Proceedings Colloque International C.N.R.S., Paris-Orsay, 1976)
[Colloques Internationaux du C.N.R.S. N$^o$ 260],
\'Editions du C.N.R.S.,
Paris, 1978, pp. 399--401.

\end{document}